\documentclass[12pt]{article} 
\usepackage{amsmath}
\usepackage{amssymb} 
\usepackage{amscd}
\usepackage{amsthm}
\usepackage[utf8]{inputenc}
\usepackage{stmaryrd}
\newtheorem{theorem}{Theorem}

\newtheorem{proposition}[theorem]{Proposition}
\newtheorem{lemma}[theorem]{Lemma}
\newtheorem{corollary}[theorem]{Corollary}

\newtheorem{definition}{Definition}

\def\dim{{\mbox{dim}}}

\def\cala{{\cal A}} 
\def\calb{{\cal B}}

\def\calh{{\cal H}}

\def\bbbone{\mbox{\rm 1\hspace {-.6em} l}}

\numberwithin{equation}{section}

\begin{document}
\enlargethispage{3cm}

\thispagestyle{empty}
\begin{center}
{\bf THE QUANTUM GROUP OF A PREREGULAR MULTILINEAR FORM}
\end{center}

\vspace{0.3cm}

\begin{center}
Julien BICHON\footnote{Laboratoire de Math\'ematiques, Universit\'e Blaise Pascal, Campus des c\'ezeaux BP 80026, F-63171 Aubi\`ere Cedex\\
Julien.Bichon@math.univ-bpclermont.fr} and
Michel DUBOIS-VIOLETTE
\footnote{Laboratoire de Physique Th\'eorique UMR 8627,
Universit\'e Paris XI,
B\^atiment 210\\ F-91 405 Orsay Cedex\\
Michel.Dubois-Violette$@$u-psud.fr} \end{center}
\vspace{0,5cm}

\begin{abstract}
We describe the universal quantum group preserving a preregular multilinear form, by means of an explicit finite presentation of the corresponding Hopf algebra. 
\end{abstract}



\section{Introduction} \label{Intro}

The quantum group of a nondegenerate bilinear form (i.e. the largest quantum group preserving the bilinear form) was introduced and constructed in \cite{mdv-lau:1990}. The aim of the paper is to generalize this construction to an appropriate class of $m$-linear forms ($m\geq 3$), the preregular multilinear forms introduced in \cite{mdv:2005,mdv:2007}. \\

Let $w$ be a preregular $m$-linear form on a vector space $V$ (see Definition \ref{PreReg}). We show the existence of a universal Hopf algebra $\mathcal H(w)$ preserving $w$ and we provide an explicit \textsl{finite} presentation by generators and relations for $\mathcal H(w)$. It should be emphasized that the existence part is something that can be considered as well-known (see the end of Section \ref{Bw})
and that our main result is that we have obtained a finite presentation for the universal object $\mathcal H(w)$.  
The condition of preregularity, introduced in \cite{mdv:2005,mdv:2007} in connection with the analysis of AS-regular algebras,  turns out to be particularly relevant in the Hopf algebra context. \\

In view of the $m=2$ case \cite{bic:2003b}
and of the construction of the quantum group $SU(m)$
as the universal compact quantum group preserving a well-chosen $m$-linear form \cite{wor:1988}, it is
natural to expect that the Hopf algebras $\mathcal H(w)$ should produce quantum analogues of $SL(m)$. We examine two natural examples that show that this is far from being true if $m\geq 3$.
 The first one is when $w$ is the ``signature'' $m$-form, for which we get a non commutative and non cosemisimple Hopf algebra, having the algebra of polynomial functions on $SL(m)$ as a quotient. This example also shows that the Hopf algebra
$\mathcal H(w)$ differs in general from those Hopf algebras $\mathcal H(w,\tilde{w})$ constructed in \cite{mdv:2007}.
The second example is when $w$ is the ``totally orthogonal'' form, and we get some of the quantum reflection groups studied in   \cite{ban-ver:2009}, whose fusion rules are noncommutative. \\

The paper is organized as follows. Section \ref{preli} is devoted to preliminaries: we fix some notation and conventions and we recall the construction of the quantum group of a nondegenerate bilinear form. In Section \ref{Prereg} we discuss preregular multilinear forms. In Section \ref{Bw} we describe the universal bialgebra preserving a multilinear form. In Section \ref{UniHw} we provide the explicit finite presentation for the universal Hopf algebra $\mathcal H(w)$ preserving a preregular multilinear form. In Section \ref{Hws} we discuss the relation between   $\mathcal H(w)$
and another class of Hopf algebras $\mathcal H(w, \tilde{w})$ introduced in \cite{mdv:2007}. Section \ref{examples} discusses some examples and the final Section 8 is our conclusion.

\section{Preliminaries}
\label{preli}

\subsection{Notations and conventions}

Throughout this article $\mathbb K$ is a (commutative) field and all vector spaces, algebras, etc. are over $\mathbb K$, the symbol $\otimes$ denotes the tensor product over $\mathbb K$ and we use the Einstein convention of summation over repeated up-down indices in the formulas. In the following $V$ is a finite-dimensional vector space with $\dim(V)=n\geq 2$ equipped with a basis $(e_\lambda)_{\lambda\in\{1,\dots, n\}}$ and we endow the dual vector space $V^\ast$ of $V$ with the dual basis $(\theta^\lambda)_{\lambda\in \{1,\dots,n\}}$ of the basis $(e_\lambda)$.  Using the finite-dimensionality of $V$ we identify $(V^{\otimes^m})^\ast$ with $(V^\ast)^{\otimes^m}$ and thus, a $m$-linear form on $V$ is the same thing as an element of $(V^\ast)^{\otimes^m}$. Given an endomorphism $Q$ of $V$, we denote by $Q^t$ its transposed that is the endomorphism of $V^\ast$ defined by $\langle \omega, QX\rangle=\langle Q^t\omega,X\rangle$ for any $\omega\in V^\ast$ and $X\in V$.

\subsection{The quantum group of a non-degenerate bilinear form}
Let $b=b_{\mu\nu}\theta^\mu\otimes \theta^\nu\in V^\ast\otimes V^\ast=(V\otimes V)^\ast$ be a nondegenerate bilinear form on $V$, ($b(e_\mu,e_\nu)= b_{\mu\nu}$). The quantum group of the nondegenerate bilinear form $b$ was defined in \cite{mdv-lau:1990} in the following manner. Denote by $b^{\mu\nu}$ the component of the inverse matrix of $(b_{\mu\nu})\in M_n(\mathbb K)$, i.e. one has
\begin{equation}
b^{\mu\lambda}b_{\lambda\nu}=b_{\nu\rho}b^{\rho\mu}=\delta^\mu_\nu
\label{binv}
\end{equation}
for $\mu,\nu\in \{1,\dots,n\}$.\\

 Note the obvious fact that $\tilde b=b^{\mu\nu} e_\mu\otimes e_\nu\in V\otimes V = (V^\ast\otimes V^\ast)^\ast$ does not depend on the basis $(e_\lambda)$ but only depends on $b\in V^\ast\otimes V^\ast=(V\otimes V)^\ast$.\\

Let $\calh(b)$ be the unital associative algebra generated by the $n^2$ elements $u^\lambda_\rho$ ($\lambda,\rho\in \{1,\dots,n\}$) with relations
\begin{equation}
b_{\mu\nu}u^\mu_\lambda u^\nu_\rho=b_{\lambda\rho}\bbbone
\label{bst}
\end{equation}
\begin{equation}
b^{\mu\nu} u^\lambda_\mu u^\rho_\nu=b^{\lambda\rho}\bbbone
\label{binst}
\end{equation}
for $\lambda,\rho\in \{1,\dots,n\}$. Then $\calh(b)$ has a unique Hopf algebra structure with coproduct $\Delta$, counit $\varepsilon$ and antipode $S$ such that
\begin{equation}
\Delta(u^\mu_\nu)=u^\mu_\lambda \otimes u^\lambda_\nu
\label{bcop}
\end{equation}
\begin{equation}
\varepsilon(u^\mu_\nu)=\delta^\mu_\nu
\label{bcou}
\end{equation}
\begin{equation}
S(u^\mu_\nu)=b^{\mu\lambda} u^\rho_\lambda b_{\rho\nu}
\label{bant}
\end{equation}
for $\mu,\nu\in \{1,\dots,n\}$. The dual object of $\calh(b)$ is {\sl the quantum group of the nondegenerate bilinear form $b$}. The analysis of the category of representations of this quantum group, that is of the category of corepresentations of $\calh(b)$, has been done in \cite{bic:2003b}.\\

Although not completely obvious with the above definition, the Hopf algebra $\calh(b)$ is characterized by the following universal property.
\begin{theorem}\label{Univb}
Let $H$ be a Hopf algebra which coacts on $V$ as
\begin{equation}
e_\lambda \mapsto e_\mu \otimes v^\mu_\lambda
\label{coH}
\end{equation}
$v^\mu_\lambda\in H$ and is such that
\begin{equation}
b_{\mu\nu} v^\mu_\lambda v^\nu_\rho=b_{\lambda\rho} \bbbone
\label{bcoH}
\end{equation}
for $\lambda, \rho\in \{1,\dots,n\}$, (where $\bbbone$ denotes the unit of $H$). Then there is a unique homomorphism of Hopf algebras $\varphi:\calh(b)\rightarrow H$ such that $\varphi(u^\lambda_\mu)=v^\lambda_\mu$ for $\lambda,\mu\in \{1,\dots,n\}$.
\end{theorem}
 
 In particular the smallest Hopf subalgebra of $H$ which contains the $v^\lambda_\mu$'s is a quotient of $\calh(b)$. For the sake of completeness we give a  proof (which will appear as a particular case of a more general result explained later, see  Section \ref{UniHw}).\\
 
 \noindent\underbar{Proof}. Relations (\ref{bcoH}) imply
 \[
 b^{\sigma\lambda} v^\mu_\lambda b_{\mu\nu} v^\nu_\rho= \delta^\sigma_\rho \bbbone
 \]
 which implies that the antipode $S$ of $H$ satisfies
 \[
 S(v^\sigma_\nu)=b^{\sigma\lambda} v^\mu_\lambda b_{\mu\nu}
 \]
 and therefore one has
 \[
 v^\sigma_\nu S(v^\nu_\rho)=v^\sigma_\nu b^{\nu\lambda} v^\mu_\lambda b_{\mu\rho}=\delta^\sigma_\rho \bbbone
 \]
from which one obtains by contraction with $b^{\rho\tau}$
\begin{equation}
b^{\nu\lambda} v^\sigma_\nu v^\tau_\lambda = b^{\sigma\tau} \bbbone
\label{coHb}
\end{equation}
for any $\sigma, \tau\in \{1,\dots, n\}$. This implies together with (\ref{bcoH}) that $u^\mu_\nu\mapsto v^\mu_\nu$ defines a homomorphism of Hopf algebras of $\calh(b)$  onto the Hopf subalgebra of $H$ generated by the $v^{\mu}_{\nu}$. $\square$\\


\section{Preregular multilinear forms}\label{Prereg}

As in the previous section $V$ is a finite-dimensional vector space with $\dim(V)=n\geq2$, $(e_\lambda)$ is a basis of $V$ with dual basis $(\theta^\lambda)$ of the dual $V^\ast$ of $V$, etc.. We shall always make the identifications $(V^{\otimes^n})^\ast=V^{\ast\otimes^n}$ which is allowed in view of the  finite dimensionality of $V$. The following definition is taken from \cite{mdv:2005} (the main motivation was Theorem 4.3 there).

\begin{definition}\label{PreReg}
Let $m$ be an integer with $m\geq2$. A $m$-linear form $w$ on $V$ is said to be {\sl preregular} iff it satisfies the following conditions $\mathrm{(i)}$ and $\mathrm{(ii)}$.\\
$\mathrm{(i)}$ If $X\in V$ is such that $w(X_1,\dots, X_{m-1},X)=0$ for any $X_1,\dots,X_{m-1}\in V$, then $X=0$.\\
$\mathrm{(ii)}$ There is an element $Q_w\in GL(V)$ such that one has
\[
w(X_1,\dots,X_{m-1},X_m)=w(Q_wX_m,X_1,\dots,X_{m-1})
\]
for any $X_1,\dots, X_m\in V$.
\end{definition}

It follows from (i) that $Q_w$ as in (ii) is unique. Condition (i) when combined with (ii) implies the stronger condition (i').\\

\noindent (i') For any $0\leq k\leq m-1$, if $X\in V$ is such that 
\[
w(X_1,\dots,X_k,X,X_{k+1},\dots,X_{m-1})=0 
\]
for any $X_1,\dots,X_{m-1}\in V$, then $X=0$.\\

Condition (i')
will be refered to as {\sl 1-site nondegeneracy} while (ii) will be refered to as {\sl twisted cyclicity} with {\sl twisting} element $Q_w$. Thus a preregular multilinear form is a multilinear form which is 1-site nondegenerate and twisted cyclic.\\

Basic examples of preregular multilinear forms are:
the signature form (see Section \ref{sign}), its $q$-analogue (defined in \cite{wor:1988}), the totally orthogonal form (see Section \ref{totort}), and the multilinear forms appearing in \S 3.3 and \S 5.3 of  \cite{mdv:2010} which correspond to representative examples of (regular) homogeneous algebras. For the latter class, it is worth noticing here that in the case of quadratic algebras, these multilinear forms give the noncommutative version of the volume element (see Proposition 10 of \cite{mdv:2007} and the comment following its proof).\\

By applying $n$ times the relations of (ii) of Definition 1 one obtains the invariance of $w$ by $Q_w$, that is
\[
w(X_1,\dots, X_m)=w(Q_wX_1,\dots,Q_wX_m)
\]
for any $X_1,\dots,X_m\in V$. Thus a twisted cyclic $w$ is invariant by its twisting element.\\

Conversely, let $w$ be an arbitrary $Q$-invariant $m$-linear form on $V$ (with $Q\in GL(V)$) then the $m$-linear form $\pi_Q(w)$ defined by 
\[
\pi_Q(w) (X_1,\dots,X_m)=\sum^m_{k=1} w(QX_k,\dots,QX_m,X_1,\dots,X_{k-1})
\]
for any $X_1,\dots,X_m\in V$ is twisted cyclic with twisting element $Q$, (in short is $Q$-cyclic). Notice that if $\mathbb K$ is of characteristic $0$, then $\frac{1}{m}\pi_Q$ is a projection of the $Q$-invariant $m$-linear forms onto the $Q$-twisted cyclic ones.\\

Condition (i) is equivalent to the existence of a $m$-linear form $\tilde w$ on $V^\ast$ such that one has
\begin{equation}
\tilde w^{\mu \lambda_1\dots\lambda_{m-1}} w_{\lambda_1\dots\lambda_{m-1}\nu}=\delta^\mu_\nu
\label{nondm}
\end{equation}
for $\mu,\nu\in \{1,\dots,m\}$, where the components of $\tilde w$ and $w$ are defined by $\tilde w^{\lambda_1\dots\lambda_m}=\tilde w(\theta^{\lambda_1},\dots,\theta^{\lambda_m})$ and $w_{\lambda_1\dots\lambda_m}=w(e_{\lambda_1},\dots, e_{\lambda_m})$ for $\lambda_k\in \{1,\dots,m\}$. Notice that in contrast with the case $m=2$,  $\tilde w$ is non unique for $m>2$. The set of all $\tilde w$ satisfying (\ref{nondm}) as above is an affine subspace of $V^{\otimes^m}$ which will be denoted by Aff($w$) and refered to as {\sl the polar of $w$}.\\

In components, condition (ii) reads
\begin{equation}
w_{\lambda_1\dots\lambda_m}=Q^\lambda_{\lambda_m} w_{\lambda\lambda_1\dots \lambda_{m-1}}
\label{Cii}
\end{equation}
for $\lambda_k\in \{1,\dots, m\}$, where the $Q^\mu_\nu$ are the components of $Q_w$ defined by
\begin{equation}
Q_w(e_\nu)=Q^\mu_\nu e_\mu
\label{compQ}
\end{equation}
for $\nu\in\{1,\dots,n\}$ in the basis $(e_\lambda)$.\\

If one combines (\ref{nondm}) with (\ref{Cii}), one gets
\[
Q^\rho_\nu \tilde w^{\mu \lambda_1\dots\lambda_{m-1}}w_{\rho\lambda_1\dots \lambda_{m-1}} = \delta^\mu_\nu
\]
that is
\begin{equation}
\tilde w^{\mu \lambda_1\dots\lambda_{m-1}}w_{\nu\lambda_1\dots \lambda_{m-1}} = (Q^{-1})^\mu_\nu
\label{exQ}
\end{equation}
for $\mu,\nu \in \{1,\dots,n\}$.\\

Let us say a few words on the case $m=2$. It is clear that a bilinear form $w=b$ satisfying Condition (i) is just a nondegenerate bilinear form and that then $\tilde w$ as in \ref{nondm} is unique and given by $\tilde w^{\mu\nu}=b^{\mu\nu}$ with the notations of Section \ref{preli}.
Hence 
a preregular bilinear form $b$ is nondegenerate. 

Conversely let $b$ be a nondegenerate bilinear form with components $b_{\mu\nu}$ in the basis $(e_\lambda)$. Then 
$b$ is twisted cyclic with
\begin{equation}
Q^\lambda_\rho=b^{\mu\lambda} b_{\mu\rho}
\label{Qb}
\end{equation}
and hence a nondegenerate bilinear form is preregular. \\

The classification of nondegenerate bilinear forms over algebraically closed fields is known: the isomorphism class of a nondegenerate bilinear form only  depends on the conjugacy class of its twisting element, see \cite{rie:1974}. On the other hand, to the best of our knowledge, no general classification result is known for $m$-linear forms if $m\geq 3$ (isomorphism of multilinear forms is called congruence in \cite{belser:2006}). The only basic observation we have is that isomorphic preregular $m$-linear forms 
have conjugate twisting elements. 


\section{The universal bialgebra $\calb(w)$}\label{Bw}

Throughout this section $w$ is a $m$-linear form on $V$ satisfying Condition (i) of Definition 1. The components of $w$ in the basis $(e_\lambda)$ are denoted by $w_{\lambda_1\dots \lambda_m}=w(e_{\lambda_1},\dots, e_{\lambda_m}$) as before.\\

Let $\calb(w)$ be the unital associative algebra generated by the $n^2$ elements $a^\mu_\nu$ ($\mu, \nu\in \{1,\dots, n\}$) with relations
\begin{equation}
w_{\lambda_1\dots\lambda_m} a^{\lambda_1}_{\mu_1}\dots a^{\lambda_m}_{\mu_m}=w_{\mu_1\dots \mu_m}\bbbone
\label{wast}
\end{equation}
for $\mu_k\in \{1,\dots,n\}$. Then $\calb(w)$ has a unique bialgebra structure with coproduct $\Delta$ and counit $\varepsilon$ such that
\begin{equation}
\Delta(a^\mu_\nu)=a^\mu_\lambda \otimes a^\lambda_\nu
\label{acop}
\end{equation}

\begin{equation}
\varepsilon(a^\mu_\nu)=\delta^\mu_\nu
\label{acou}
\end{equation}
for $\mu, \nu\in\{1,\dots,n\}$. Let $\tilde w$ ($\in V^{\otimes^m}$) be, as in Section 2, a solution of (\ref{nondm}). One has in view of (\ref{wast})
\begin{equation}
\tilde w^{\mu \lambda_1\dots\lambda_{m-1}} a^{\rho_1}_{\lambda_1} \dots a^{\rho_{m-1}}_{\lambda_{m-1}} w_{\rho_1\dots \rho_{m-1}\sigma}  a^\sigma_\nu=\delta^\mu_\nu
\label{leftinva}
\end{equation}
for $\mu,\nu\in \{1,\dots,n\}$, which means that in $M_n(\calb(w))$ the matrix $(a^\mu_\nu)$ has $(\tilde w^{\mu\lambda_1\dots \lambda_{m-1}} a^{\rho_1}_{\lambda_1}\dots a^{\rho_{m-1}}_{\lambda_{m-1}}w_{\rho_1\dots \rho_{m-1}\nu})$ as left inverse. 

\begin{theorem}\label{UnivBw}
Let $B$ be a bialgebra which coacts on $V$ as
\[
e_\lambda\mapsto e_\mu \otimes v^\mu_\lambda
\]
$v^\mu_\nu \in B$ and is such that
\[
w_{\rho_1\dots \rho_m} v^{\rho_1}_{\lambda_1} \dots v^{\rho_m}_{\lambda_m}=w_{\lambda_1\dots \lambda_m} \bbbone
\]
for $\lambda_k\in \{1,\dots,n\}$, (where $\bbbone$ is the unit of $B$). Then there is a unique homomorphism of bialgebras $\varphi:\calb(w)\rightarrow B$ such that 
\[
\varphi(a^\mu_\nu)=v^\mu_\nu\>\>\>\>\>\>   \forall  \mu,\nu\in \{1,\dots, n\}.
\]
If furthermore $B$ is a Hopf algebra with antipode $S$, then one has 
\[
S(v^\mu_\nu)=\tilde w^{\mu \lambda_1\dots\lambda_{m-1}} v^{\rho_1}_{\lambda_1} \dots v^{\rho_{m-1}}_{\lambda_{m-1}} w_{\rho_1\dots \rho_{m-1} \nu} 
\]
for $\mu, \nu \in \{1,\dots, n\}$, with $\tilde w$ as above.
\end{theorem}

\noindent \underbar{Proof}. The first part is straightforward and classical in a more general context while the second part follows from
\[
\begin{array}{ll}
& \tilde w^{\mu \lambda_1\dots\lambda_{m-1}} v^{\rho_1}_{\lambda_1} \dots v^{\rho_{m-1}}_{\lambda_{m-1}} w_{\rho_1\dots \rho_{m-1}\sigma}  v^\sigma_\tau S(v^\tau_\nu)\\
\\
& = S(v^\mu_\nu)=\tilde w^{\mu \lambda_1\dots\lambda_{m-1}} v^{\rho_1}_{\lambda_1} \dots v^{\rho_{m-1}}_{\lambda_{m-1}} w_{\rho_1\dots \rho_{m-1}\nu}
\end{array}
\]
for $\mu,\nu\in \{1,\dots,n\}$ which follows from (\ref{leftinva}) and from the definition of the antipode. $\square$ \\

Hence by Manin's construction of the Hopf envelope of a bialgebra \cite{man:1988},
there exists a universal Hopf algebra coacting on a multilinear form $w$, which is finitely
generated as an algebra since $w$ satisfies condition (i) in Definition \ref{PreReg}
(by the last identies of the previous theorem). We have not been able to provide a finite presentation of this algebra without using the full definition of preregularity.

\section{The universal Hopf algebra $\calh(w)$}\label{UniHw}

Throughout this section $w$ is a given preregular multilinear form on $V$ with twisting element $Q_w$ simply denoted by $Q\in GL(V)$.\\

The components of $w, Q$ and $Q^{-1}$ in the basis $(e_\lambda)$ of $V$ (with dual basis $(\theta^\lambda)$ of $V^\ast$) are denoted by $w_{\lambda_1\dots\lambda_m},Q^\mu_\nu$ and $(Q^{-1})^\mu_\nu$ with
\[
w_{\lambda_1\dots \lambda_m} = w(e_{\lambda_1},\dots,e_{\lambda_m})
\]
\[
Q(e_\nu)=Q^\mu_\nu e_\mu,\ \ \ \text{or}\ \ \ \  Q^\mu_\nu=\langle \theta^\mu,Q(e_\nu)\rangle
\]
\[
Q^{-1}(e_\nu)=(Q^{-1})^\mu_\nu e_\mu,\ \ \ \text{or}\ \ \ \  (Q^{-1})^\mu_\nu=\langle \theta^\mu,Q^{-1}(e_\nu)\rangle
\]
for $\lambda_k,\mu,\nu\in \{1,\dots,n\}$.\\

Let $\calh(w)$ be the unital associative algebra generated by the $2n^2$ generators $u^\mu_\nu, s^\mu_\nu$ with relations
\begin{equation}
u^\mu_\lambda s^\lambda_\nu =\delta^\mu_\nu \bbbone
\label{us}
\end{equation}
\begin{equation}
Q^\lambda_\nu u^\rho_\lambda(Q^{-1})^\sigma_\rho s^\mu_\sigma =\delta^\mu_\nu\bbbone
\label{tus}
\end{equation}
\begin{equation}
w_{\lambda_1\dots \lambda_m} u^{\lambda_1}_{\mu_1}\dots u^{\lambda_m}_{\mu_m}=w_{\mu_1\dots \mu_m} \bbbone
\label{invw}
\end{equation}
for $\mu,\nu,\mu_k\in \{1,\dots,n\}$. 

\begin{proposition}\label{HwA}
One has in $\calh(w)$
\begin{equation}
w_{\mu_1\dots \mu_m} s^{\mu_m}_{\nu_m}\dots s^{\mu_1}_{\nu_1} = w_{\nu_1\dots \nu_m}\bbbone
\label{sinw}
\end{equation}
\begin{equation}
s^\mu_\rho u^\rho_\nu=\delta^\mu_\nu \bbbone
\label{su}
\end{equation}
\begin{equation}
s^\lambda_\nu Q^\tau_\lambda u^\rho_\tau (Q^{-1})^\mu_\rho=\delta^\mu_\nu \bbbone
\label{tsu}
\end{equation}
for $\nu_k,\mu,\nu \in \{1,\dots,n\}$. Furthermore, for any solution $\tilde w$ of $(\ref{nondm})$ i.e. $\forall \tilde w\in \mathrm{Aff}(w)$ one has
\begin{equation}
s^\mu_\nu=\tilde w^{\mu \lambda_1\dots\lambda_{m-1}}u^{\rho_1}_{\lambda_1}\dots u^{\rho_{m-1}}_{\lambda_{m-1}}w_{\rho_1\dots\rho_{m-1}\nu}
\label{Rsu}
\end{equation}
\begin{equation}
Q^\tau_\nu u^\rho_\tau (Q^{-1})^\mu_\rho=\tilde w^{\mu \lambda_1\dots\lambda_{m-1}}s^{\rho_{m-1}}_{\lambda_{m-1}}\dots s^{\rho_1}_{\lambda_1}w_{\rho_1\dots\rho_{m-1}\nu}
\label{Rus}
\end{equation}
for $\mu,\nu\in \{1,\dots,n\}$.
\end{proposition}

\noindent \underbar{Proof}. Relations (\ref{sinw}) are obtained from (\ref{invw}) by contraction on the right by $s^{\mu_m}_{\nu_m}\dots s^{\mu_1}_{\nu_1}$ and by using (\ref{us}). From (\ref{nondm}) and (\ref{invw}) it follows that 
\[
(\tilde w^{\mu \lambda_1\dots\lambda_{m-1}}u_{\lambda_1}^{\rho_1}\dots u^{\rho_{m-1}}_{\lambda_{m-1}} w_{\rho_1\dots \rho_{m-1}\sigma}) u^\sigma_\nu=\delta^\mu_\nu \bbbone
\]
 which implies (\ref{Rsu}) in view of (\ref{us}) and then (\ref{su}) follows (see also Theorem \ref{UnivBw}). Similarily, from (\ref{nondm}) and (\ref{sinw}) it follows that 
\[
s^\sigma_\nu(\tilde w^{\mu \lambda_1\dots\lambda_{m-1}} s^{\rho_{m-1}}_{\lambda_{m-1}} s^{\rho_1}_{\lambda_1}w_{\rho_1\dots \rho_{m-1}\sigma})=\delta^\mu_\nu \bbbone
\]
which implies (\ref{Rus}) in view of (\ref{tus}) and then (\ref{tsu}) follows. $\square$

\begin{proposition} \label{HopfHw}
The algebra $\calh(w)$ has a unique structure of Hopf algebra with coproduct $\Delta$, counit $\varepsilon$ and antipode $S$ such that
\begin{equation}
\Delta(u^\mu_\nu) = u^\mu_\lambda\otimes u^\lambda_\nu
\label{copu}
\end{equation}
\begin{equation}
\Delta(s^\mu_\nu) = s^\lambda_\nu\otimes s_\lambda^\mu
\label{cops}
\end{equation}
\begin{equation}
\varepsilon(u^\mu_\nu) = \varepsilon(s^\mu_\nu)=\delta^\mu_\nu
\label{cou}
\end{equation}
\begin{equation}
S(u^\mu_\nu)=s^\mu_\nu
\label{antu}
\end{equation}
\begin{equation}
S(s^\mu_\nu)=(Q^{-1})^\mu_\rho u^\rho_\lambda Q^\lambda_\nu 
\label{ants}
\end{equation}
for $\mu,\nu\in \{1,\dots,n\}$.
\end{proposition}

\noindent \underbar{Proof}. This is a direct verification using the relations from the previous proposition. $\square$ \\




We did not use Condition (ii) of Definition \ref{PreReg} (twisted cyclicity) to construct the Hopf algebra $\mathcal H(w)$ ($Q$ could be any invertible matrix). The twisted cyclicity condition is fully used to prove the following universal property of $\mathcal H(w)$.


\begin{theorem}\label{Uniw}
Let $w$ be preregular $m$-linear form and let $H$ be a Hopf algebra which coacts on $V$ as
\[
e_\lambda \mapsto e_\mu \otimes v^\mu_\lambda
\]
$v^\mu_\lambda\in H$ and is such that
\[
w_{\mu_1\dots \mu_m}v^{\mu_1}_{\lambda_1}\dots v^{\mu_m}_{\lambda_m}=w_{\lambda_1\dots\lambda_m}\bbbone
\]
for $\lambda_k\in \{1,\dots,n\}$, (where $\bbbone$ denotes the unit of $H$). Then there is a unique homomorphism of Hopf algebras $\varphi:\calh(w)\rightarrow H$ such that $\varphi(u^\lambda_\mu)=v^\lambda_\mu$ for $\lambda,\mu\in \{1,\dots,n\}$.
\end{theorem}
This theorem is the announced generalization of Theorem \ref{Univb}. It follows that if $w =b$ is bilinear, then $\mathcal H(b)$ as in Section \ref{preli} is isomorphic with $\mathcal H(w)$ as above. \\

\noindent\underbar{Proof}. All what one has to verify is that the relations (\ref{us}), (\ref{tus}) and (\ref{invw}) are satisfied with $u^\mu_\nu$ replaced by $v^\mu_\nu$ and $s^\mu_\nu=S(u^\mu_\nu)$ replaced by $S(v^\mu_\nu)$, ($S$ being the antipode). Relation (\ref{invw}) is satisfied by definition and relation (\ref{us}) is a consequence of the axioms for Hopf algebras. Thus it remains to show that
\begin{equation}
Q^\lambda_\nu v^\rho_\lambda(Q^{-1})^\sigma_\rho S(v^\mu_\sigma)=\delta^\mu_\nu \bbbone
\label{tvs}
\end{equation}
is satisfied in $H$. On the other hand, Theorem \ref{UnivBw} implies that 
\[
S(v^\mu_\nu)=\tilde w^{\mu\lambda_1\dots \lambda_{m-1}} v^{\rho_1}_{\lambda_1} \dots v^{\rho_{m-1}}_{\lambda_{m-1}}  w_{\rho_1\dots \rho_{m-1}\nu}
\]
with $\tilde w$ solution of (\ref{nondm}). Therefore one has
\[
\begin{array}{ll}
Q^\lambda_\nu v^\rho_\lambda(Q^{-1})^\sigma_\rho S(v^\mu_\sigma) & =  Q^\lambda_\nu v^\rho_\lambda (Q^{-1})^\sigma_\rho \tilde w^{\mu\lambda_1\dots \lambda_{m-1}} v^{\rho_1}_{\lambda_1} \dots v^{\rho_{m-1}}_{\lambda_{m-1}} w_{\rho_1\dots \rho_{m-1}\sigma}\\
&=\tilde w^{\mu\lambda_1\dots \lambda_{m-1}}Q^\lambda_\nu v^\rho_\lambda v^{\rho_1}_{\lambda_1} \dots v^{\rho_{m-1}}_{\lambda_{m-1}} w_{\rho \rho_1\dots \rho_{m-1}}\\
&= \tilde w^{\mu\lambda_1\dots \lambda_{m-1}} Q^\lambda_\nu w_{\lambda\lambda_1\dots \lambda_{m-1}}\bbbone\\
&=\tilde w^{\mu\lambda_1\dots\lambda_{m-1}}w_{\lambda_1\dots \lambda_{m-1}\nu}\bbbone = \delta^\mu_\nu\bbbone
\end{array}
\]
where the twisted cyclicity condition is used in the second equality. This proves (\ref{tvs}). $\square$ \\

It follows in particular from Theorem \ref{Uniw} that the algebra of polynomial functions on ${\rm Aut}(w)$ is a quotient of $\mathcal H(w)$. An element of ${\rm Aut}(w)$ has to commute with the twisting element $Q$, so if $Q$ is diagonalizable with distinct eigenvalues, ${\rm Aut}(w)$ is very small, being a subgroup of the torus $\mathbb (K^*)^m$. This still allows $\mathcal H(w)$ to be a large algebra (and hence $w$ has a large quantum symmetry group), as the example of quantum $SL(2)$ shows.

\section{The Hopf algebras $\calh(w,\tilde w)$, $\tilde w\in \mathrm{Aff}(w)$}\label{Hws}

In this section, where $w$ is a preregular $m$-linear form on $V$ with twisting element $Q\in GL(V)$ and $\tilde w$ is an element of $\mathrm{Aff}(w)$ (that is an element of $V^{\otimes^m}$ satisfying (\ref{nondm})),
we recall the definition of the Hopf algebra $\mathcal H(w,\tilde{w})$, constructed in \cite{mdv:2007}, and relate it to $\mathcal H(w)$.

Let $\calh(w,\tilde w)$ be the unital associative algebra generated by the $n^2$ elements $v^\mu_\nu$ ($\mu, \nu\in \{1,\dots,n\}$) with relations
\begin{equation}
w_{\lambda_1\dots\lambda_m} v^{\lambda_1}_{\mu_1}\dots v^{\lambda_m}_{\mu_m}=w_{\mu_1\dots \mu_m}\bbbone
\label{wvst}
\end{equation}
and
\begin{equation}
\tilde w^{\lambda_1\dots\lambda_m} v^{\mu_1}_{\lambda_1}\dots v^{\mu_m}_{\lambda_m}=\tilde w^{\mu_1\dots \mu_m}\bbbone
\label{tildewv}
\end{equation}
for $\mu_k\in \{1,\dots,n\}$. There is a unique structure of Hopf algebra on $\calh(w,\tilde w)$ with coproduct $\Delta$, counit $\varepsilon$ and antipode $S$ such that
\begin{equation}
\Delta v^\mu_\nu=v^\mu_\lambda \otimes v^\lambda_\nu
\label{vcop}
\end{equation}
\begin{equation}
\varepsilon(v^\mu_\nu)=\delta^\mu_\nu
\label{vcou}
\end{equation}
and
\begin{equation}
S(v^\mu_\nu)=\tilde w^{\mu\lambda_1\dots\lambda_{m-1}}v^{\rho_1}_{\lambda_1}\dots v^{\rho_{m-1}}_{\lambda_{m-1}}w_{\rho_1\dots \rho_{m-1}\nu}
\label{vant}
\end{equation}
for $\mu, \nu\in \{1,\dots, n\}$.\\

By Theorem \ref{Uniw}  there exists a unique (surjective) homomorphism of Hopf algebras
\[
\varphi:\calh(w)\rightarrow \calh(w,\tilde w)
\]
such that $\varphi(u^\mu_\nu)=v^\mu_\nu$ for $\mu,\nu\in\{1,\dots,n\}$.
We will see in the next section that this is not an isomorphism in general.



\section{Examples}
\label{examples}

In this section, in order to illustrate our construction, we examine two examples. 
We begin  with a general useful lemma, whose proof is straightforward combining
Relations (\ref{sinw}) and (\ref{su}).

\begin{lemma}\label{relations}
Assume that $m\geq 3$. We have, in $\mathcal H(w)$,
$$w_{\lambda \rho \mu_3 \ldots \mu_{m}} s^{\mu_m}_{\nu_m}
\ldots s^{\mu_{3}}_{\nu_{3}}=w_{\nu_1 \ldots \nu_{m}}
u^{\nu_1}_{\lambda}u^{\nu_{2}}_\rho$$
for any $\lambda$, $\rho$,  $\nu_3, \ldots , \nu_{m} \in \{1, \ldots , n\}$.
\end{lemma}

\subsection{Example 1: the signature form}\label{sign}

Let $w=\varepsilon$ be the signature form on $V=\mathbb K^m$ (the volume form), i.e.
$\varepsilon_{\lambda_1\ldots \lambda_m}=0$ if two indices are equal, and 
$\varepsilon_{\lambda_1\ldots \lambda_m}$ is the signature of the corresponding permutation otherwise. This is a preregular $m$-linear form, with $Q=I_m$.

\begin{proposition}
 The Hopf algebra $\mathcal H(\varepsilon)$ is not commutative if $m\geq 3$. 
\end{proposition}
 
\noindent\underbar{Proof}. It is easy to see that it is enough to prove the result at $m=3$. Let $\mathcal A$ be the free algebra on two generators $x$ and $y$. It is straightforward to check the existence of an algebra map $\mathcal H(w)\rightarrow \mathcal A$
$$
\begin{pmatrix}
 u^1_1 & u^1_2 & u^1_3 \\
u^2_1 & u^2_2 & u^2_3 \\
u^3_1 & u^3_2 & u^3_3 
\end{pmatrix}
\mapsto
\begin{pmatrix}
 1 & x & y \\
0 & 1 & 0 \\
0 & 0 & 1 
\end{pmatrix}, \ 
\begin{pmatrix}
 s^1_1 & s^1_2 & s^1_3 \\
s^2_1 & s^2_2 & s^2_3 \\
s^3_1 & s^3_2 & s^3_3 
\end{pmatrix}
\mapsto
\begin{pmatrix}
 1 & -x & -y \\
0 & 1 & 0 \\
0 & 0 & 1 
\end{pmatrix}
$$
Thus $\mathcal H(\varepsilon)$ is not commutative. $\square$ \\

Using Lemma \ref{relations} we see that the following relations hold in $\mathcal H(\varepsilon)$:
$$u^{\lambda}_\nu u^\mu_\nu = u^\mu_\nu u^\lambda_\nu, \ [u^\lambda_\nu, u^\mu_\rho]
=[u^\mu_\nu, u^\lambda_\rho]$$
for any $\lambda, \mu, \nu, \rho$. These are the relations defining the Manin matrices studied in \cite{cfr:2009}.

Let $\tilde{\varepsilon}$ be given by $\tilde{\varepsilon}_{\lambda_1\ldots \lambda_m}
=\frac{1}{m}\varepsilon_{\lambda_1\ldots \lambda_m}$ (assuming here that $m \not=0$ in $\mathbb K$). We have $\tilde{\varepsilon} \in {\rm Aff}(\varepsilon)$.

\begin{proposition}
Assume that $m \geq 3$. Then the canonical Hopf algebra map $\mathcal H(\varepsilon) \rightarrow \mathcal H(\varepsilon, \tilde{\varepsilon})$
is not injective, and $\mathcal H(\varepsilon)$ is not cosemisimple.
\end{proposition}

\noindent\underbar{Proof}.
One can  show that $\mathcal H(\varepsilon, \tilde{\varepsilon})$
is commutative, as in \cite{wor:1988,ros:1990}. It thus follows that the canonical Hopf algebra map $\mathcal H(\varepsilon) \rightarrow \mathcal H(\varepsilon, \tilde{\varepsilon})$
is not injective. 

Assume that $\mathcal H(\varepsilon)$ is cosemisimple. Then there exists $w' \in {\rm Aff}(\varepsilon)$ satisfying (\ref{tildewv}). But $\tilde{\varepsilon}$
is, up to scalar, the only element in ${\rm Aff}(\varepsilon)$ for which (\ref{tildewv}) can hold in 
$\mathcal H(\varepsilon, \tilde{\varepsilon})$ (by some well-known facts on the representation theory of $SL(m)$). Thus, up to a scalar, $w'= \tilde{\varepsilon}$, and relations (\ref{tildewv}) hold for $\tilde{\varepsilon}$ in $\mathcal H(\varepsilon)$. Thus the canonical map $\mathcal H(\varepsilon) \rightarrow \mathcal H(\varepsilon, \tilde{\varepsilon})$ is an isomorphism: contradiction.
$\square$

\subsection{Example 2: the totally orthogonal form}\label{totort}

We now consider the totally orthogonal $m$-linear form $\theta$  on $V=\mathbb K^n$:
$\theta_{\lambda_1\ldots \lambda_m}=1$ if all indices are equal and $0$ otherwise.
This is a preregular $m$-linear form with $Q=I_n$.

\begin{proposition}
 For $m\geq 3$, the Hopf algebra $\mathcal H(\theta)$ is isomorphic to the algebra
presented by generators $a^{\mu}_\lambda$, $\mu, \lambda \in \{1, \ldots , n\}$, and relations
$$a^\mu_\lambda a^\mu_\nu=0=a^\lambda_\mu a^\nu_\mu, \ \forall \mu, \ \forall \lambda \not=\nu$$
$$\sum_{\lambda=1}^n (a^\mu_\lambda)^m=1= \sum_{\lambda=1}^n (a^\lambda_\mu)^m, \ \forall \mu $$
Thus $\mathcal H(\theta)$ is isomorphic with the Hopf algebra denoted $A_h^m(n)$ \cite{ban-ver:2009}.
\end{proposition}

\noindent\underbar{Proof}.
Let $\mathcal A$ denote the algebra having the presentation given in the statement of the proposition. It is a direct verification to check that there exists a unique Hopf algebra structure on $\mathcal A$ such that 
$$
\Delta(a^\mu_\nu) = a^\mu_\lambda\otimes a^\lambda_\nu, \ 
\varepsilon(a^\mu_\nu) =\delta^\mu_\nu, \ 
S(a^\mu_\nu)=(a^\nu_\mu)^{m-1}$$
and that $\mathcal A$ satisfies the assumption in Theorem \ref{Uniw}. Hence there exists a surjective Hopf algebra map $\mathcal H(\theta) \rightarrow \mathcal A$, $u^\mu_\lambda \mapsto a^\mu_\lambda$.

To construct the inverse isomorphism, we have to check that the defining relations in $\mathcal A$ are satisfied in $\mathcal H(\theta)$.
The second family of relations holds in $\mathcal H(\theta)$ by Relations (\ref{invw}). 
We have, by Lemma \ref{relations}, $u^\mu_\lambda u^\mu_\nu=0$ for $\lambda \not=\nu$.
Using the antipode we get  $s^\mu_\lambda s^\mu_\nu=0$ for $\lambda \not=\nu$.
Moreover, using Relations (\ref{invw}), we see that $s^\lambda_\mu=(u_\lambda^\mu)^{m-1}$, hence $s_\mu^\lambda s_\mu^\nu=0$ for $\lambda \not=\nu$, and using the antipode, we see that the first family of relations holds in $\mathcal H(\theta)$. We get a surjective algebra map $\mathcal A \rightarrow \mathcal H(\theta)$, $a^\mu_\lambda \mapsto u^\mu_\lambda$, and we are done. $\square$ \\

The Hopf algebra $\mathcal H(\theta)\simeq A_h^m(n)$ (which also can be described by a free wreath product
operation \cite{bic:2004}) is, when $\mathbb K=\mathbb C$, cosemimple and its fusion rules are non-abelian, see \cite{ban-ver:2009}. It would be interesting to characterize the preregular multilinear forms for which the diagrammatic
 techniques in \cite{ban-ver:2009} can be extended.

\section{Conclusion : The quantum group of a preregular multilinear form}

The Hopf algebra $\calh(w,\tilde w)$ which was considered in \cite{mdv:2007} in connection with its coaction on the algebras $\cala(w,N)
$  defined by the ``twisted potential" $w$ (for  $m\geq N\geq 2$) has the defect to be dependent of the auxiliary variable $\tilde w\in \text{Aff}\ (w)$, (excepted for $m=2$ for which $\text{card}(\text{Aff}\ (w))=1)$.\\

 The Hopf algebra $\calh(w)$ coacts as well on the algebras $\cala(w,N)$ defined in Section 5 of \cite{mdv:2007}. This latter class of algebras includes all AS-regular algebras which are $N$-Koszul as shown in \cite{mdv:2007}, Theorem 11 (see also in \cite{mdv:2005}).\\

 It is clear from Theorem \ref{Uniw} that the right Hopf algebra which only depends on $w$ is the universal Hopf algebra $\calh(w)$ defined in Section \ref{UniHw}. Therefore we define {\sl the quantum group of the preregular multilinear form $w$} to be the dual object of $\calh(w)$. The {\sl representations} of this quantum group are the corepresentations of the Hopf algebra $\calh(w)$ which plays the role of the Hopf algebra of ``representative functions" on it.\\
 
 In fact the Hopf algebras $\calh(w,\tilde w)$ for $\tilde w\in \text{Aff}(w)$ all are quotients of $\calh(w)$, which means that the corresponding quantum groups are quantum subgroups of the above quantum group of the preregular multilinear form $w$. In the case $m=2$, that is when $w$ is a nondegenerate bilinear form, then $\tilde w$ is unique and $\calh(w,\tilde w)$ coincides with $\calh(w)$ as explained in Section \ref{preli} (Theorem \ref{Univb}). Thus for $m=2$, the canonical projection
 \[
 \pi:\calh(w)\rightarrow \calh(w,\tilde w)
 \]
 is an isomorphism. However for $m\geq 3$, Example 1 of the previous section
shows that this not always true, but Example 2 shows it can be true. 
Those examples also show that the uniform description of the 
corepresentation category of $\mathcal H(w)$ for the case $m=2$ in \cite{bic:2003b} will no longer
hold when $m\geq 3$. \\

The Hopf algebras $\calh(w)$ and $\calh(w,\tilde w)$ above belong to a particular class of Hopf algebras named cosovereign Hopf algebras in \cite{bic:2001} where they were considered and studied in a categorical perspective. Examples of Hopf algebras of this type are those arising from compact quantum groups \cite{wor:1987}. The class of cosovereign Hopf algebras is also exactly the class of Hopf algebras considered in \cite{ac-mos:1998} (see also in \cite{cra:2002}) and the analysis of \cite{ac-mos:1998} was later extended to a larger class of Hopf algebras, which mixed the notion of cosovereign Hopf algebra with the dual notion of sovereign Hopf algebra \cite{bic:2001}, in \cite{ac-mos:1999} and \cite{ac-mos:2000}. These classes of Hopf algebras generalize the class of Hopf algebras with antipode satisfying $S^2=I$ for which cyclic cohomology and the generalization of the Weil homomorphism have tractable descriptions.

\bibliographystyle{plain}

\end{document}